\documentclass[11pt]{article}
\usepackage{amsmath,amssymb}

\newtheorem{theo}{Theorem}[section]
\newtheorem{prop}[theo]{Proposition}
\newtheorem{lemma}[theo]{Lemma}

\newtheorem{coro}[theo]{Corollary}
\newtheorem{conj}[theo]{Conjecture}
\newcommand\npf{\mbox{ }\hfill\sqr\vskip6pt}
\def\sqr{$\vcenter{\hrule height.2mm
\hbox{\vrule width.2mm height2mm\kern2mm
\vrule width.2mm}\hrule height.2mm}$}

\begin{document}

\title{Edge coloring complete uniform hypergraphs with many components}

\author{
Yair Caro
\thanks{Department of Mathematics, University of Haifa at Oranim, Tivon 36006, Israel. e--mail: yairc@macam98.ac.il}
\and
Raphael Yuster
\thanks{Department of Mathematics, University of Haifa at Oranim, Tivon 36006, Israel. e--mail: raphy@macam98.ac.il}
}

\date{}

\maketitle

\begin{abstract}
Let $H$ be a hypergraph. For a $k$-edge coloring $c : E(H) \rightarrow
\{1,\ldots,k\}$ let $f(H,c)$ be the number of components in the
subhypergraph induced by the color class with the least number of
components. Let $f_k(H)$ be the maximum possible value of $f(H,c)$
ranging over all $k$-edge colorings of $H$. If $H$ is the complete graph $K_n$ then, trivially,
$f_1(K_n)=f_2(K_n)=1$. In this paper we prove that for $n \geq 6$,
$f_3(K_n)=\lfloor n/6 \rfloor+1$ and supply close upper and lower
bounds for $f_k(K_n)$ in case $k \geq 4$. Several results concerning the value of
$f_k(K_n^r)$, where $K_n^r$ is the complete $r$-uniform hypergraph on $n$ vertices, are also
established.
\end{abstract}

\setcounter{page}{1}

\section{Introduction}
All graphs and hypergraphs considered here are finite, unordered and simple. For
standard terminology the reader is referred to \cite{We}. Let
$H$ be a hypergraph. For a $k$-edge coloring $c : E(H) \rightarrow
\{1,\ldots,k\}$ let $f(H,c)$ be the number of components in the
subhypergraph induced by the color class with the least number of
components. Isolated vertices are not considered as components in
a subhypergraph induced by edges. To avoid trivialities we always
assume $k \leq e(H)$. Let $f_k(H)$ be the maximum possible value
of $f(H,c)$ ranging over all $k$-edge colorings of $H$. Trivially,
if $H$ has $n$ vertices then $1 \leq f_k(H) \leq \lfloor n/r
\rfloor$ where $r$ is the minimum cardinality of an edge of $H$,
and $f_1(H)=c(H)$ where $c(H)$ is the number of
components of $H$ that are not isolated vertices.

In case $H$ is the complete $r$-uniform hypergraph $K_n^r$ it is not difficult to
prove (see last part of Theorem \ref{t3}) that $f_k(K_n^r)=1$ for $k \leq r$. In fact,
this is a generalization of the graph theoretic case, $r=2$, where $f_2(K_n)=1$
is merely the well-known fact that either a graph or its complement is connected
(recently several extension of this
elementary fact were proved \cite{BiDiVo,CaRo}
and our paper is in part inspired by this simple theorem).
For fixed $k \geq r+1$ it is
not difficult to show that $f_k(K_n^r)$ is linear in $n$. However,
determining the exact value is a nontrivial task.

This paper contains several results concerning the parameter $f_k(K_n^r)$.
In the graph-theoretic case, we completely
settle the case $k=3$ and the cases $k=n-1$ and $k=n$.
For other fixed values of $k$ we supply close upper and
lower bounds that are also valid for the hypergraph case.

Our main results are summarized in the following
theorems. For simplicity we use the notation $f(n,k,r)$ instead of
$f_k(K_n^r)$, and $f(n,k)$ for the graph-theoretic case $r=2$.
The first theorem on $f(n,k)$ is an exact result dealing with the lower end of the scale of colors, namely $k=3$.
\begin{theo}
\label{t1}
$f(4,3)=f(5,3)=2$. Otherwise $f(n,3)=\lfloor n/6 \rfloor +1$.
\end{theo}
The next theorem is an exact result dealing with values of $k$ in the upper end of
the scale.
\begin{theo}
\label{t2}
$~$
\begin{enumerate}
\item
$f(n,n-1)=\lfloor n/2 \rfloor$.
\item
$f(n,n)=\lfloor (n-1)/2 \rfloor$.
\item
If $k \geq n-1$ and $k$ divides ${n \choose 2}$ then $f(n,k)=n(n-1)/(2k)$.
\item
Suppose $t \geq 1$ and $rt$ divides $n$, then for $k=t(n-1)!/((r-1)!(n-r)!)$,
$f(n,k,r)=n/(rt)$.
\item
If $k \geq {n \choose r} - {{n-r} \choose r}$
then $f(n,k,r)=\lfloor {n \choose r}/k \rfloor$.
\end{enumerate}
\end{theo}

Our next theorem supplies close upper and lower bounds for all
fixed values of $k$. Before we state the theorem we need a
few definitions. Given an edge-coloring $c$ of $K_n^r$, let $z(n,r,c,s)$ denote the
fraction of the vertices incident with at least one edge whose color
is $s$. Let $z(n,r,c)$ denote the maximum value of $z(n,r,c,s)$
taken over all colors appearing in $c$.
For $k \geq 1$, let $z_k(n,r)$ denote the minimum possible value
of $z(n,r,c)$ taken over all colorings that use at most $k$ colors.
Finally, let $z_{k,r}$ denote the infimum of $z_k(n,r)$ taken over all $n \geq r$.
For $r=2$ denote $z_k=z_{k,2}$. The following theorem relates $z_{k,r}$ with the function $f(n,k,r)$.
\begin{theo}
\label{t3}
Let $k \geq r+1$. Then,
$$
n\left(\frac{1}{r}-\frac{1}{rk^{1/r}}\right)(1+o_n(1)) \geq
f(n,k,r) \geq
n\left(\frac{1}{r}-\frac{z_{k,r}}{r}\right)(1+o_n(1)).
$$
If $k \leq r$ then $z_{k,r}=1$ and $f(n,k,r)=1$.
\end{theo}
We conjecture that the lower bound in the last theorem is the correct one:
\begin{conj}
\label{c1}
Let $k \geq r+1$. Then,
$$
f(n,k,r) = n\left(\frac{1}{r}-\frac{z_{k,r}}{r}\right)(1+o_n(1)).
$$
\end{conj}
In section 2 we analyze the parameters $z_{k,r}$ and $z_k$.
Infinitely many values of $z_{k,r}$ are known, and the values of
infinitely many others are open problems.
In particular, we determine (with varying difficulty of proofs depending on $k$)
the following specific values: $z_{r+1,r}=r/(r+1)$, $z_4=3/5$, $z_5=5/9$, $z_6=1/2$, $z_7=3/7$, $z_{12}=1/3$,
$z_{6,3}=2/3$, $z_{14,3}=1/2$. More generally,
$z_{p^2+p+1}=(p+1)/(p^2+p+1)$ whenever $p$ is a prime power,
and $z_{p^2+p}=1/p$ whenever $p$ is a prime power.
$k=8$ is the smallest number for which $z_k$ is still unknown.
It is also not difficult to bound $z_k$ from below.
In fact, we show $z_k \geq 1/\lceil \sqrt{k+1/4}-1/2\rceil$.
Since, by definition, $z_k$ is monotone decreasing and since
for every integer $s$ there is a prime (moreover a prime power) between $s$ and $s+O(s^{2/3})$
\cite{Hu} we have that $z_k=1/\sqrt{k}+O(1/k)$.
Together with Theorem \ref{t3} we have, in the graph-theoretic case:
\begin{coro}
Let $k \geq 3$ and let $p$ be the largest prime power satisfying
$p^2+p+1 \leq k$ then
$$
n\left(\frac{1}{2}-\frac{1}{2\sqrt{k}}\right)(1+o_n(1)) \geq
f(n,k) \geq
n\left(\frac{1}{2}-\frac{p+1}{2(p^2+p+1)}\right)(1+o_n(1)).
$$
\end{coro}
The last corollary, together with the argument concerning density of primes
show that the upper and lower bounds in the last corollary are very close.
Table \ref{T1} summarizes the best upper and lower bounds that we currently have
for $\limsup f(n,k)/n$ and $\liminf f(n,k)/n$ respectively, for some small values
of $k$. For each specific $k$ (except $k=3$), the upper bound follows from Theorem
\ref{t3} and in all cases, the lower bound follows
from the best known upper bound for $z_k$ (all these upper bounds are consequences of constructions
that appear in Section 2). We note that it
may be possible to prove Conjecture \ref{c1} without
determining the precise value of $z_k$ or $z_{k,r}$ for all $k$. On the other hand, although we know
that, say, $z_4=3/5$, and $z_{6,3}=2/3$, proving that $\lim f(n,4)/n = 1/5$ or $\lim f(n,6,3)/n=1/9$
is still an open problem.
Currently, Conjecture \ref{c1} is open, in the graph theoretic case, for all $k \geq 4$, and for $r \geq 3$
it is open for all $k \geq r+1$.

\begin{table}[ht]
\begin{center}
\begin{tabular}{|c|c|c|c|c|}
\hline
$k$ & $\limsup f(n,k)/k \leq$ & $\liminf f(n,k)/k \geq$ & $z_k \leq$ & $z_k \geq$ \\ \hline
$3$ & $\frac{1}{6}$ & $\frac{1}{6}$ & $\frac{2}{3}$ & $\frac{2}{3}$\\
$4$ & $\frac{1}{4}$ & $\frac{1}{5}$ & $\frac{3}{5}$ & $\frac{3}{5}$ \\
$5$ & $0.277$ & $0.222$ & $\frac{5}{9}$ & $\frac{5}{9}$ \\
$6$ & $0.296$ & $0.25$ & $\frac{1}{2}$ & $\frac{1}{2}$\\
$7$ & $0.311$ & $0.285$ & $\frac{3}{7}$ & $\frac{3}{7}$ \\
$8$ & $0.324$ & $0.285$ & $\frac{3}{7}$ & $\frac{3}{8}$ \\
$9$ & $\frac{1}{3}$ & $\frac{3}{10}$ & $\frac{2}{5}$ & $\frac{1}{3}$ \\
$10$ & $0.342$ & $0.3$ & $\frac{2}{5}$ & $\frac{1}{3}$ \\
$11$ & $0.35$ & $0.318$ & $\frac{4}{11}$ & $\frac{1}{3}$ \\
$12$ & $0.356$ & $0.333$ & $\frac{1}{3}$ & $\frac{1}{3}$ \\
$13$ & $0.362$ & $0.346$ & $\frac{4}{13}$ & $\frac{4}{13}$ \\ \hline
\end{tabular}
\end{center}
\caption{\label{T1}Asymptotic upper and lower bounds for small values of $k$}
\end{table}

In Section 3 we prove the theorems and also consider the bipartite analog,
namely, the parameter $f_k(K_{n,n})$.

\section{Localized edge-coloring of complete hypergraphs}

\subsection{Lower bounds for $z_{k,r}$}

We consider first the case $k \leq r$. We show that in this case $z_{k,r}=1$.
Namely, there must be a color that is incident with all vertices. In fact, we will show
something somewhat stronger. There is always a color such that the subgraph of
$K_{n,r}$ induced by the edges with this color is connected and spanning. In other
words:
\begin{lemma}
\label{l20}
If $ k \leq r$ then $f(n,k,r)=1$.
\end{lemma}
{\bf Proof:}\,
Clearly it suffices to prove $f(n,r,r)=1$. We shall prove this
by induction on $r$. For $r=2$ we trivially have $f(n,2)=f(n,2,2)=1$ since
either a graph or its complement is connected. Assume the theorem holds for $r-1$,
and we prove it for $r$. The proof for $r$ also proceeds by induction on $n$.
If $n=r+1$ then exactly one color appears in two (intersecting) edges and hence this color must be incident with $r+1=n$
vertices and we trivially have $f(r+1,r,r)=1$. Assume it hold for $n-1$ and we prove it for $n$.
Fix a vertex $v$ and consider the subhypergraph $K_{n-1}^r$ obtained by removing $v$.
By the induction hypothesis there is a color $c$ such that the subhypergraph of $K_{n-1}^r$ induced
by this color is spanning and connected. If $c$ also appears in some edge of $K_n^r$ that contains
$v$ we are done. Otherwise, consider the hypergraph $K_{n-1}^{r-1}$ obtained by removing $v$ and
coloring each $(r-1)$-set $f$ with the color of the edge $f \cup {v}$ in $K_n^r$.
This coloring does not use the color $c$ and hence is an $r-1$ coloring of $K_{n-1}^{r-1}$.
By the induction hypothesis it has a color $c'$ such that the subhypergraph induced by this color
is connected and spanning. By the definition of the coloring of $K_{n-1}^{r-1}$,
the subhypergraph of $K_n^r$ induced by color $c'$ must also be connected and spanning. \npf.

The following lemma supplies a lower bound for $z_{k,r}$ in case $k \geq r+1$.
\begin{lemma}
\label{l21}
Let $k \geq r+1$, and let $d \geq 2$ be a positive integer.
Then,
$$
z_{k,r} \geq \max \left\{
\min \left\{ \left(\frac{d}{k}\right)^{\frac{1}{r-1}}~,~\frac{1}{d-1}\right\}~,~
\min \left\{ \frac{d}{k}~,~z_{d-1,r-1}\right\}
\right\}.
$$
\end{lemma}
{\bf Proof:}\,
Let $d \geq 2$ be a positive integer.
Consider a coloring of $K_n^r$ with at most $k$ colors.
For $i=1,\ldots,k$, let $s_i$ denote the number of $(r-1)$-sets that are contained in an edge
whose color is $i$. For an $(r-1)$-set $U$, let $b_U$
denote the number of distinct colors that appear in edges that contain $U$.
Clearly, $\sum_{U \in {{[n]} \choose {r-1}}}b_U=s_1+\cdots+s_k$.
If $b_U \geq d$ for all $U \in {{[n]} \choose {r-1}}$
then $s_1+\cdots+s_k \geq {n \choose {r-1}}d$ and hence at least some $s_i$ satisfies $s_i \geq {n \choose {r-1}}d/k$.
Thus, if $x$ is the number of vertices appearing in an edge colored $i$, we must have
${x \choose {r-1}} \geq s_i \geq {n \choose {r-1}}d/k$.
This implies that
$$
\left(\frac{x}{n}\right)^{r-1} \geq \frac{x(x-1) \cdots (x-r+2)}{n(n-1) \cdots (n-r+2)} \geq \frac{d}{k}
$$
and hence $x \geq n (d/k)^{1/(r-1)}$.
Thus, the fraction of vertices incident with color $i$ is at least $(d/k)^{1/(r-1)}$.
Otherwise, there exists an $(r-1)$-subset $U$ such that $b_U \leq d-1$. Assume, without loss of
generality, that only the colors $1,\ldots,d-1$ appear in edges that contain $U$.
Then each vertex of $K_n^r$ appears in an edge whose color is one of the colors
$1,\ldots,d-1$. Thus, for some color $i$
the fraction of vertices appearing in an edge colored $i$ is at least $1/(d-1)$.
We have shown that $z_{k,r} \geq \min\{(d/k)^{1/(r-1)}~,~1/(d-1)\}$.

\noindent
Next, we show that $z_{k,r} \geq \min\{d/k, z_{d-1,r-1}\}$.
For completeness, define $z_{t,1}=1/t$ (in a $t$-coloring of $n$ singletons
at least $\lceil n/t \rceil$ elements obtain the same color, and this bound is realized).
If there is a vertex $v$ of $K_n^r$ that is incident only with $t \leq d-1$
colors then consider a coloring of $K_{n-1}^{r-1}$ obtained by removing
vertex $v$ and coloring each $(r-1)$-set $U$ with the color of the edge $U \cup \{v\}$ in
$K_n^r$. This defines a $t$-coloring of $K_{n-1}^{r-1}$ and
hence, by definition of $z_{t,r-1}$, some color is incident with at least $(n-1)z_{t,r-1}$ vertices.
Since this color is also incident with $v$ and since $z_{d-1,r-1} \leq z_{t,r-1}$ we have
found a color incident with a fraction of at least $z_{d-1,r-1}$ vertices of $K_n^r$.
Otherwise, each vertex of $K_n^r$ is incident with at  least $d$ colors, and hence
there is a color that is incident with at least $nd/k$ vertices. In any case we have shown
$z_{k,r} \geq \min\{d/k, z_{d-1,r-1}\}$.
\npf

In the case $r=2$, by choosing $d=\lceil \sqrt{k+1/4}+1/2\rceil$ and using the fact that
for this choice of $d$, $d/k \geq 1/(d-1)$, we get
\begin{coro}
\label{c21}
For all $k \geq 3$, $z_k \geq 1/\lceil \sqrt{k+1/4}-1/2\rceil$.
\end{coro}
Notice that for $r \geq 3$ it is best to take $d \approx k^{1/r}$ in Lemma \ref{l21}, except for
small values of $k$ where the lower bound of $\min\{d/k, z_{d-1,r-1}\}$ does better.

\subsection{Upper bounds and precise values of $z_{k,r}$}
Upper bounds for $z_k$ and $z_{k,r}$ are demonstrated by construction.
Consider a trivial coloring of $K_n^r$ with ${n \choose r}$ colors,
each edge colored with a unique color. The fraction of vertices
incident with each color is trivially $r/n$. This shows, in particular,  that
\begin{coro}
\label{c22}
$z_{{n \choose r},r} \leq r/n$ ~, ~$z_{n(n-1)/2} \leq 2/n$~,~ $z_{r+1,r} \leq r/(r+1)$. \npf
\end{coro}
Corollary \ref{c22}, together with Lemma \ref{l21} yield the following:
\begin{coro}
\label{c23}
$z_3=2/3$ ~,~ $z_{r+1,r} = r/(r+1)$~,~$z_6=1/2$.
\end{coro}
{\bf Proof:}\,
Upper bounds follow from Corollary \ref{c22}. $z_{r+1,r} \geq r/(r+1)$
follows from Lemma \ref{l21} by taking $k=r+1$ and $d=r$ (and recalling that
$z_{r-1,r-1}=1$ by Lemma \ref{l20}). $z_3$ is a special case of $z_{r+1,r}$.
$z_6 \geq 1/2$ by taking $r=2$, $k=6$ and $d=3$ in Lemma \ref{l21}. \npf

In many cases we can find non-trivial constructions that match the lower bound that follows from
Lemma \ref{l21}. For example, the smallest nontrivial Steiner Triple System shows that
$K_7$ can be decomposed into $7$ triangles. In other words, there is a coloring
of $K_7$ with $7$ colors such that each color induces a triangle. We therefore get $z_7 \leq 3/7$.
On the other hand, applying
Lemma \ref{l21} with $r=2$, $k=7$ and $d=3$ gives $z_7 \geq 3/7$. Thus, $z_7=3/7$.
More generally we can prove the following:
\begin{prop}
\label{p20}
Let $p$ be a prime power. Then,
$z_{p^2+p+1}=(p+1)/(p^2+p+1)$ and $z_{p^2+p} = 1/p$.
\end{prop}
{\bf Proof:}\,
Whenever $p$ is a prime power there exists a projective plane $PG(2,p)$.
This projective plane corresponds
to the existence of a $2-(p^2+p+1,p+1,1)$ design (see, e.g., \cite{CoDi}),
which, in turn, corresponds to the fact that
$K_{p^2+p+1}$ decomposes into $p^2+p+1$ copies of $K_{p+1}$.
Hence, we have that $z_{p^2+p+1} \leq (p+1)/(p^2+p+1)$
and using Lemma \ref{l21} with $r=2$ and $d=p+1$ we get $z_{p^2+p+1} = (p+1)/(p^2+p+1)$.
Similarly, when $p$ is a prime power there exists an affine plane $AG(2,p)$.
This affine plane
 corresponds to the existence of a $2-(p^2,p,1)$ design (\cite{CoDi}),
which, in turn, corresponds to the fact that
$K_{p^2}$ decomposes into $p^2+p$ copies of $K_p$
and using Lemma \ref{l21} with $r=2$, $d=p+1$ we get $z_{p^2+p} = 1/p$. \npf

\noindent
Notice that Proposition \ref{p20} gives, in particular, $z_{13}=4/13$ and $z_{12}=1/3$.

Hanani has shown the existence of $3-(n,4,1)$ designs for every even $n$ not divisible by 6
(cf. \cite{CoDi}). In other words, $K_4^3$ decomposes $K_n^3$ for $n=2,4 \bmod 6$.
Hence if $k=n(n-1)(n-2)/24$ and $n=2,4 \bmod 6$ then $z(k,3) \leq  4/n$.
In particular $z(14,3) \leq 1/2$. This upper bound has a matching lower bound that
follows from Lemma \ref{l21} by taking $r=3$, $k=14$ and $d=7$ and recalling
that $z_{6,2}=z_6=1/2$. Thus, we obtain the sporadic value $z_{14,3}=1/2$.

In all previous constructions, all color classes induced the same clique (or hyperclique) size.
Constructions using non-isomorphic color classes (and even color classes that are not cliques)
are also very useful. In fact, sometimes using non-isomorphic color classes is provably an optimal strategy.
Consider the case $k=4$. Color $K_5$ with four colors as follows:
For $i=1,2,3$, color $i$ appears in the edges $(i,4)$ and $(i,5)$
and color $4$ appears in the edges $(1,2),(1,3),(2,3)$. The edge $(4,5)$ is colored
arbitrarily by one of the colors $1,2$ or $3$. Notice that each color is incident with precisely
three vertices. Thus, $z_4 \leq 3/5$. Notice that we cannot match this
upper bound with Lemma \ref{l21} so in order to prove that this is an optimal strategy
we need to explicitly prove:
\begin{prop}
\label{p1}
$z_4 \geq 3/5$.
\end{prop}
{\bf Proof:}\,
We need to show that in any coloring of a complete graph
$K_n$ with at most four colors, at least one color is incident with at least $3n/5$ vertices.
Consider a coloring of $K_n$ with the colors $1,2,3,4$.
We use the same notations as in the proof of Lemma \ref{l21} (recalling that in case $r=2$ the
sets $U$ are singletons).
If $b_v \geq 3$ for all $v=1,\ldots,n$ then at least one color has $s_i \geq 3n/4 > 3n/5$.
Assume, therefore, that there exist vertices with $b_v=2$. (If there exists a vertex with $b_v=1$
then the unique color $i$ incident with $v$ has $s_i=n$.) Let $X \subset \{1,\ldots,n\}$
be the subset of vertices with $b_v=2$. Each $v \in X$ is associated with the unique pair of colors
incident with $v$. Notice that if $v$ is associated with $(i,j)$ then no $v' \in X$ is associated with
$(k,l)$ where $\{k,l\} \cap \{i,j\} = \emptyset$. Thus, there are only at most three types of associations.
It follows that there are two (not necessarily mutually exclusive) cases.
Either there exists a color $i$ such that each $v \in X$ is incident with $i$,
or there exists a color $j$ such that no $v \in X$ is incident with $j$. Consider the first case.
If $|X| \geq 3n/5$ then $s_i \geq 3n/5$ and we are done. Otherwise, $s_1+s_2+s_3+s_4 \geq 3(n-|X|)+2|X|=3n-|X|\geq 12n/5$.
Thus, some $s_i$ has $s_i \geq 3n/5$. Consider the second case. Then
$s_1+s_2+s_3+s_4-s_j \geq 2n$. Hence, some $s_i$ (i $\neq j$) has $s_i \geq 2n/3$. \npf

The coloring of $K_5$ with four colors described above yielding $z_4 \leq 3/5$ is somewhat
''reducible''. Indeed, if we allow weights on the vertices of $K_n$ then the following coloring
of $K_4$ with four colors is more efficient: Color the edge $(i,4)$ with color $i$ for $i=1,2,3$.
Color the triangle $(1,2,3)$ with color 4. Assign the weight $2/5$ to vertex 4, and the weight $1/5$
to each of the vertices $1,2,3$. Thus, the sum of weights of the vertices incident with color $i$
is $3/5$ for all $i=1,2,3,4$, as expected. Indeed, any $c$-coloring of a weighted $K_n$ (with rational weights)
can be transformed to a non-weighted $c$-coloring of a larger $K_{n'}$ where each vertex of $K_n$
is ``blown-up'' proportionally to its weigh (an edge of $K_{n'}$ connecting two vertices corresponding to the same
blown-up vertex of $K_n$ can be colored with any arbitrary color that is incident with that vertex of $K_n$).

The proof for the case $k=5$ is more complicated. For the upper bound,
color $K_9$ as follows: Two colors induce each a copy of $K_5$. The two copies of $K_5$ share one vertex.
The remaining 16 yet uncolored edges form a $K_{4,4}$. $K_{4,4}$ can be decomposed into two $K_{2,3}$ and
one $K_{1,4}$. Thus, we have a 5-coloring of $K_9$ where each color is incident with 5 vertices.
This shows that $z_5 \leq 5/9$. Again, we cannot match this upper bound with Lemma \ref{l21}.
We therefore show:
\begin{prop}
\label{p2}
$z_5 \geq 5/9$.
\end{prop}
{\bf Proof:}\,
We need to show that in any coloring of a complete graph
$K_n$ with at most five colors, at least one color is incident with at least $5n/9$ vertices.
Consider a coloring of $K_n$ with the colors $1,2,3,4,5$.
We use the same notations as in the proof of Lemma \ref{l21}.
If some vertex has $b_v=1$ we are done.
Assume, therefore all vertices have $b_v \geq 2$. If at least $7n/9$ vertices have $b_v > 2$
then $s_1+s_2+s_3+s_4+s_5 \geq 3n-2n/9$ which implies that some $i$ has $s_i \geq 5n/9$ as required.
Thus, let $X \subset \{1,\ldots,n\}$ be the subset of vertices with $b_v=2$, and we may assume
$|X| \geq 2n/9$. As in Proposition \ref{p1}, each $v \in X$ is associated with the unique pair of colors
incident with $v$, and if $v$ is associated with $(i,j)$ then no $v' \in X$ is associated with
$(k,l)$ where $\{k,l\} \cap \{i,j\} = \emptyset$. Thus,
either there exists a color $i$ such that each $v \in X$ is incident with $i$,
or there exist three colors $i,j,k$ such that only the associations $(i,j),(i,k),(j,k)$
are valid associations for the vertices of $A$.
Consider the second case. In this case, each $u \notin X$ is incident with at least two of the
colors $\{i,j,k\}$. To see this, notice that if $u$ is not incident with, say, both $i$ and $j$ then the
edge connecting $u$ to $v \in A$ where $v$ is associated with $(i,j)$, cannot be colored.
We therefore have $s_i+s_j+s_k \geq 2n$ so at least one of the colors is incident with  at least $2n/3 > 5n/9$
vertices.
Now consider the first case. Since color $i$ is incident with all vertices of $X$ we may assume $|X| \leq 5n/9$,
otherwise we are done.
Let $s$ be the number of additional colors, other than $i$,  incident with vertices of $X$.
Without loss of generality, assume $i=1$ and assume all other colors incident with vertices of $X$ are $2,\ldots,s+1$.
Clearly $1 \leq s \leq 4$. Let $Y$ be the subset of vertices not in $X$ and which are incident with
color 1. Thus, color 1 is incident with precisely $|Y|+|X|$ vertices, and colors $2,\ldots,s+1$ are all incident with
at least all vertices of $v(K_n) \setminus (X \cup Y)$, namely to at least $n-|X|-|Y|$ vertices.
Assume first that $s \geq 2$. Each vertex of $Y$ is incident with at least three colors, one of which is
color 1, so it is also incident with at least $s-2$ of the colors $2,\ldots,s+1$.
It follows that $s_2+\ldots+s_{s+1} \geq s(n-|X|-|Y|)+(s-2)|Y|+|X|$.
Thus, some $s_i$ is incident with at least $n-|X|-2|Y|/s+|X|/s$ vertices.
It suffices to show that $\max\{|Y|+|X|~, ~n-|X|-2|Y|/s+|X|/s\}\geq 5n/9$.
Since the maximum is minimized when both terms are equal, and this happens
when $|Y|=ns/(s+2)-|X|(2s-1)/(s+2)$ it suffices to show that
$ns/(s+2)-|X|(s-3)/(s+2) \geq 5n/9$. For $s=2$ this holds since $|X| \geq 2n/9$ and thus $n/2+|X|/4 \geq 5n/9$.
For $s=3$ this holds since $3n/5 > 5n/9$. For $s=4$ this holds since $|X| \leq 5n/9$ and thus
$2n/3-|X|/6 > 5n/9$. Finally, assume $s=1$. In this case colors 1 and 2 are incident with all
vertices of $K_n$ and therefore $s_1+s_2 \geq n+|X|$. Thus, one of them is incident with at least
$n/2+|X|/2 > 5n/9$ since $|X| \geq 2n/9$. \npf

The reader may notice the significant added complexity to the proof of Proposition \ref{p2} as opposed to
Proposition \ref{p1}. It is plausible that with an increasing amount of effort one may determine $z_k$ or $z_{k,r}$
for every specific $k$ and $r$ with an appropriate ``ad-hoc'' proof.

A non-symmetric construction in the hypergraph ($r=3$) case that yields an exact result is the following:
\begin{prop}
\label{p3}
$z_{6,3}=2/3$.
\end{prop}
{\bf Proof:}\,
The lower bound follows from Lemma \ref{l21} by taking $r=3$, $k=6$ and $d=4$ and
recalling that $z_{3,2}=z_3=2/3$. For the upper bound, color $K_6^3$ with six colors
as follows:
Color 1 appears in $(1,2,3),(1,2,4),(1,3,4),(2,3,4)$.
Color 2 appears in $(1,2,5),(1,2,6),(1,5,6),(2,5,6)$.
Color 3 appears in $(1,3,5),(1,3,6),(3,5,6)$.
Color 4 appears in $(2,3,5),(2,4,5),(3,4,5)$.
Color 5 appears in $(1,4,5),(1,4,6),(4,5,6)$.
Color 6 appears in $(2,3,6),(2,4,6),(3,4,6)$.
Every color is incident with four vertices hence $z_{6,3} \leq 4/6=2/3$. \npf

The upper and lower bounds for $z_k$ in case $k=8,9,10,11$ appearing in Table \ref{T1}
are obtained as follows: For $k=8$, $z_8 \leq z_7 =3/7$, and $z_8 \geq 3/8$ by selecting $d=3$
in Lemma \ref{l21}. For $k=9$, $z_9 \leq 2/5$ since $K_4^-$ (the graph obtained from $K_4$ be deleting an edge)
decomposes $K_{10}$ into $9$ copies (see \cite{CoDi}).
$z_9 \geq 1/3$ by selecting $d=4$ in Lemma \ref{l21}.
$z_{10} \leq z_9 \leq 2/5$. $z_{10} \geq 1/3$ by selecting $d=4$ in Lemma \ref{l21}.
$z_{11} \leq 4/11$ since $K_4^-$ decomposes $K_{11}$ into $11$ copies (see \cite{CoDi}).
$z_{11} \geq 1/3$ by selecting $d=4$ in Lemma \ref{l21}.

Finally, infinitely many nontrivial upper bounds in the case $r=3$ are obtained
by using M\"obius designs. For $q$ a prime power there exist the M\"obius designs $3-(q^s,q+1,1)$
(cf. \cite{CoDi}). Hence, $K_{q+1}^3$ decomposes $K_{q^s+1}^3$ and hence, for $s=2$
we have that for $k=(q^2+1)q$, $z(k,3) \leq (q+1)/(q^2+1)$. In particular,
$z(10,3) \leq 3/5$. This can be compared with the lower bound of
$z(10,3) \geq 5/9$ that follows from Lemma \ref{l21} for $r=3$, $k=10$ and $d=6$,
and using Proposition \ref{p2} which states that $z_{5,2}=z_5=5/9$.

\section{Proof of the main results}
We start with the lower bound in theorem \ref{t3} as its proof will also be used
to prove one direction in Theorem \ref{t1}.

\subsection{A lower bound for $f(n,k,r)$}
We shall require the following lemma, whose proof is an immediate consequence of
a theorem of Baranyai \cite{Ba}.
\begin{lemma}
\label{l30}
Let $r \geq 2$ and $t \geq 1$ be positive integers. There exists $n_0=n_0(r,t)$ such that
for all $n \geq n_0$, $K_n^r$ contains $t$ disjoint maximum matchings
(a maximum matching is a set of $\lfloor n/r \rfloor$ independent edges). \npf
\end{lemma}
It is well-known that $n_0(2,t)=t+1$, since the complete graph on $t+1$
vertices decomposes into $t$ perfect matchings in case $t$ is odd, and $t+1>t$ maximum matchings
in case $t$ is even. For $r \geq 3$ the theorem of Baranyai states that if $n$ is a multiple of $r$ then
$K_n^r$ decomposes into perfect matchings, and hence precisely $t=(n-1)!/((r-1)!(n-r)!)$ disjoint perfect
matchings. Thus, $n_0(r,t) \approx (t(r-1)!)^{1/(r-1)}$.

We need to prove that $f(n,k,r) \geq n\left(\frac{1}{r}-\frac{z_{k,r}}{r}\right)(1+o_n(1))$.
Let $\epsilon > 0$. Let $t=t(k,r,\epsilon)$ be the smallest integer such that
$K_t^r$ has a coloring $C$ using at most $k$ colors, such that each color is
incident with at most $t(z_{k,r}+\epsilon)$ vertices.
Notice that in some cases it is possible to
have $\epsilon=0$. In Section 2 we have shown, e.g., that for $r=2$ and $k=3,4,5,6,7,12,13$ we can take $\epsilon=0$
(in other words, $z_k$ is {\em realized} in these cases).
As an example, for $r=2$ and $k=3$ we can pick $\epsilon=0$ and $t=3$.
\begin{lemma}
\label{l31}
Suppose $n$ satisfies
$$
\left\lfloor \frac{n}{t} \right\rfloor \geq n_0\left(r ~,~k-\left\lceil \frac{1}{z_{k,r}+\epsilon} \right\rceil\right),
$$
then $f(n,k,r) \geq \lfloor \frac{n}{rt} \rfloor \lceil t(1-z_{k,r}-\epsilon)\rceil+1$.
\end{lemma}
{\bf Proof:}\,
Assume the vertices of $K_t^r$ are $\{1,\ldots,t\}$
and assume the colors of $C$ are $\{1,\ldots,k\}$.
Let $C(e)$ denote the color of edge $e$ (notice that it is possible
that $C$ is not onto $\{1,\ldots,k\}$ as $c$ may use less than
$k$ colors).

Partition the vertices of $K_n^r$ into $t$ equitable parts $A_1,\ldots,A_t$.
The cardinality of each part is either $\lceil n/t \rceil$ or $\lfloor n/t \rfloor$ and by the
assumption in the statement of the lemma, the parts are nonempty.

We now show how to color each edge of $K_n^r$. Let $f$ be an arbitrary edge,
and let $U = \{i ~:~A_i \cap f \neq \emptyset\}$. Clearly $1 \leq |U| \leq r$.
Consider first the case $|U| > 1$. In this case, let $e \in K_t^r$ be an edge such that
$U \subset e$. In this case, we color $f$ with the color $C(e)$.
Consider next the case $|U|=1$. In this case $f$ is completely within some part $A_i$.
Hence, it remains to show how to color edges that are completely within some $A_i$.
Let $C_i \subset \{1,\ldots,k\}$ be the subset of colors not appearing in any edge of
$K_t^r$ that contains the vertex $i$. We claim that $|C_i| \leq 
k-\lceil \frac{1}{z_{k,r}+\epsilon} \rceil$.
To see this, notice that if there were less than $\frac{1}{z_{k,r}+\epsilon}$
colors incident with $i$ then at least one color would have been incident with more than
$t(z_{k,r}+\epsilon)$ vertices, contradicting the assumption.
We claim that we can find in $A_i$ a set of $|C_i|$ disjoint maximum matchings. Indeed this
follows from Lemma \ref{l30}, and by the assumption in the current lemma that states that
$|A_i| \geq \lfloor n/t \rfloor \geq n_0(r,k-\lceil \frac{1}{z_{k,r}+\epsilon} \rceil)$. Fixing $|C_i|$
disjoint maximum matchings in $A_i$ we now color each maximum matching with a distinct color of $C_i$.
It remains to show how to color the edges completely within $A_i$ that do not
belong to any of the selected maximum matchings. In this case, we can color them with
any color of $C \setminus C_i$ (does not matter which one).

Now consider any color $c \in C$. If $c$ is not used at all in the coloring of $K_t^r$ then the subhypergraph
of $K_n^r$ induced by the edges colored $c$ is composed of isolated edges, exactly $\lfloor |A_i|/r \rfloor$
isolated edges from each $A_i$. Therefore, this subhypergraph has at least $t\lfloor n/(rt) \rfloor$ components.
If $c$ is used as a color of at least one edge of $K_t^r$ then the subhypergraph induced
by the edges colored $c$ in $K_n^r$ has some large components consisting of all the $s$ sets $A_i$ such that
the vertex $i$ of $K_t^r$ is incident with an edge colored $c$,
and  isolated edges, exactly $\lfloor |A_i|/r \rfloor$
isolated edges from each of the $t-s$ sets $A_i$ where vertex $i$ of $K_t^r$ is not incident
with any edge colored $c$.
Thus, this subhypergraph has at least $(t-s)\lfloor n/(rt) \rfloor+1$ components.
Since $s \leq (z_{k,r}+\epsilon)t$ our construction shows that $f(n,k,r) \geq \lfloor \frac{n}{rt} \rfloor \lceil t(1-z_{k,r}-\epsilon)\rceil+1$.
\npf

\subsection{Proof of Theorem \ref{t1}}
Trivially $f(3,3)=1$. $f(4,3)=2$ as seen by coloring $K_4$ with three colors each inducing a perfect
matching. $f(5,3)=2$ as seen by coloring $K_5$ with a red triangle and a red edge vertex-disjoint with the triangle,
and coloring the remaining 6 edges with three blue edges forming a $K_{1,2}$ and a $K_{1,1}$
vertex disjoint with each other, and three green edges that are now also forced to induce a
$K_{1,2}$ and a $K_{1,1}$ vertex disjoint with each other. Thus, we assume $n \geq 6$.

For the lower bound we use lemma \ref{l31}.
Recall that $z_3=2/3$ and $K_3$ realizes $z_3$ with $\epsilon=0$ by giving each edge
of $K_3$ a distinct color.
The condition in Lemma \ref{l31} is satisfied for all $n \geq 6$ since $n_0(2,1)=2$, and
hence we get $f(n,3) \geq \lfloor n/6 \rfloor + 1$.

It remains to show the upper bound.
Consider a coloring of $K_n$ with the three colors $1,2,3$.
Let $A_i$ be the subset of vertices in a maximum cardinality component induced by color $i$.
Put $|A_i|=a_i$ and assume $a_1 \geq a_2 \geq a_3$.
Put $D=A_1 \cap A_2$ and $d=|D|$. Put $M=V \setminus (A_1 \cup A_2)$ where
$V$ is the set of vertices of $K_n$, and $m=|M|$.
If $a_1 \geq 2n/3$ then we are done since the subgraph induced by color $1$
has at most $\lfloor 1+(n-2n/3)/2 \rfloor =1+\lfloor n/6 \rfloor$ components.
Hence we assume $a_i < 2n/3$ for $i=1,2,3$.
Consider the following cases
\begin{enumerate}
\item
$d=0$. In this case every edge between $A_1$ and $A_2$ is colored 3.
Thus, $a_3 \geq a_1+a_2$, a contradiction since the $a_i$ cannot be zero.
\item
$d >0$, $m=0$, and $d=a_1$ or $d=a_2$ (or both).
In this case one of the colors 1 or 2 (or both) induces a connected spanning subgraph
of $K_n$.
\item
$d > 0$, $m=0$, $d < a_1$, $d < a_2$.
In this case all edges between $A_1 \setminus D$ and $A_2 \setminus D$ are colored 3.
Thus, $2n/3 > a_3 \geq (a_1-d)+(a_2-d)$. On the other hand $(a_1-d)+(a_2-d)+d=n$.
Thus, $d > n/3$. It follows that $a_1+a_2=n+d> 4n/3$. Hence, $a_1 > 2n/3$, a contradiction.
\item
$d > 0$, $m > 0$, $a_1 > d$ and $a_2 > d$. As in the previous case all edges between $A_1 \setminus D$ and $A_2 \setminus D$
are colored 3. Also all edges between $M$ and $D$ are colored 3.
Thus, the subgraph induced by $3$ has at most two components and $2 \leq \lfloor n/6 \rfloor +1$.
\item
$d > 0$, $m > 0$, $a_2=d$. In this case all edges between $M$ and $D$ are colored 3.
Thus, $a_3 \geq d+m > d=a_2$, a contradiction to the assumption $a_3 \leq a_2$.
 \npf
\end{enumerate}

\subsection{Proof of Theorem \ref{t2}}
\begin{enumerate}
\item
If $n$ is even then $K_n$ decomposes into $n-1$ perfect matchings.
Thus, $f(n,n-1) \geq n/2$. If $n$ is odd then $K_n$ decomposes
into $(n-1)/2$ Hamiltonian cycles. Each cycle further decomposes into
a matching of size $(n-1)/2$ and the remaining $(n+1)/2$ edges form
a subgraph with $(n-1)/2$ components (one component is a path with three
vertices and the others are independent edges). Thus, $f(n,n-1) \geq (n-1)/2=\lfloor n/2 \rfloor$.
In both the even and odd cases we always have $f(n,n-1) \leq \lfloor
{n \choose 2}/(n-1) \rfloor=\lfloor n/2 \rfloor$.
Thus, $f(n,n-1)=\lfloor n/2 \rfloor$.
\item
If $n$ is odd then the chromatic index of $K_n$ is $n$.
Thus, $K_n$ decomposes into $n$ matchings. This forces each matching to be
of size $(n-1)/2$. Thus, $f(n,n) \geq (n-1)/2$.
If $n$ is even then $n+1$ is odd and we have $f(n+1,n)=n/2$
by the previous case. Deleting the additional vertex causes each induced subgraph of
a color to loose at most one component. Thus, $f(n,n)\geq n/2-1$.
In both cases we always have $f(n,n) \geq \lfloor (n-1) /2 \rfloor$.
Since, trivially, $f(n,n) \leq \lfloor {n \choose 2}/n \rfloor=\lfloor (n-1) /2 \rfloor$
we have shown $f(n,n)=\lfloor (n-1)/2 \rfloor$.
\item
Assume $k \geq n-1$ and assume $t=n(n-1)/(2k)$ is an integer.
Thus, $2t \leq n$ and it is well known that in this case $tK_2$ decomposes $K_n$ \cite{He,FoFu}.
Thus, $f(n,k) \geq t$. The other direction follows from the trivial fact that
$f(n,k) \leq {n \choose 2}/k=t$.
\item
Suppose $t \geq 1$ and $rt$ divides $n$. In particular, $r$ divides $n$, and by the result of
Baranyai \cite{Ba}, $K_n^r$ decomposes into $(n-1)!/((r-1)!(n-r)!)$ perfect matchings.
Each perfect matching has cardinality $n/r$ and is further decomposed into $t$ matchings,
each consisting of $n/(rt)$ edges. Hence, for $k=t(n-1)!/((r-1)!(n-r)!)$ we have
$f(n,k,r) \geq n/(rt)$. Trivially, $f(n,k,r) \leq {n \choose r}/k=n/(rt)$. \npf
\item
Consider the line graph $L$ of $K_n^r$. $L$ has $N={n \choose r}$ vertices and
is regular of degree $s={n \choose r} - {{n-r} \choose r}-1$.
By the theorem of Hajnal and Szemer\'edi \cite{HaSz}, for all $k' > s$, $L$ has a $k'$ equipartite coloring,
namely, a vertex-coloring with $k'$ colors where each color class has either $\lceil N/k' \rceil$
or $\lfloor N/k' \rfloor$ elements. Using $k=k'$ and translating this back to the original hypergraph
$K_n^r$, we have a $k$-edge coloring of $K_n^r$ such that each color induces a matching of
cardinality at least $\lfloor {n \choose r}/k \rfloor$. Thus, $f(n,k,r) \geq \lfloor {n \choose r}/k \rfloor$.
The other direction is trivial. \npf
\end{enumerate}

\subsection{Proof of Theorem \ref{t3}}
In case $k \geq r+1$, the lower bound of Theorem \ref{t3} is shown in Lemma \ref{l31}.
For the upper bound, notice that in a $k$-edge coloring of $K_n^r$ at least
one color appears in at least ${n \choose r}/k$ edges.
Fix such a color, and let $s_1,\ldots,s_t$ denote the cardinalities of the
components in the subhypergraph induced by this color.
Clearly, $s_i \geq r$ for $i=1,\ldots,t$ and
\begin{equation}
\label{e1}
{{n-r(t-1)} \choose r} +(t-1) \geq {{s_1+\cdots+s_t-r(t-1)} \choose r} +(t-1) \geq 
\sum_{i=1}^t {{s_i} \choose r} \geq {n \choose r} \cdot \frac{1}{k}.
\end{equation}
It follows that $n-r(t-1) \geq (1+o_n(1))n/k^{1/r}$ and hence
$t \leq (n/r-n/(rk^{1/r}))(1+o_n(1))$.
In case $k \leq r$, $f(n,k,r)=1$ is shown in Lemma \ref{l20}. \npf

It is not difficult to show that in case $r=2$, the negation of the error term in the upper bound of Theorem \ref{t3}
is at least $\Theta(1/k^{3/2})$, by considering the maximal cardinality components induced by
each color, and showing that either there is a color with a ``huge'' component, or otherwise there must
be (at least) two such maximal components (belonging to
two distinct colors) that intersect in $\Theta(n)$ vertices, and thus at least one of these maximal components
is far from being a clique (misses $\Theta(n^2)$ edges from being a clique), and hence the convexity
argument in inequality (\ref{e1}) cannot be exploited to its extreme.
However, this improvement is negligible. Even for $k=4$ we could only improve
the constant from $0.25$ to (a little less than) $0.243$. That is,
$\limsup f(n.4)/4 \leq 0.243$.

\subsection{The bipartite analog}
Let $C$ be a $k$-edge coloring of $K_r$ where color $i$ is incident with
$s_i$ vertices. We show how to construct a coloring $C'$ of $K_{r,r}$
where color $i$ is incident with $2s_i$ vertices. Assume the vertices
of one partite class are labeled $\{a_1,\ldots,a_r\}$ and
the vertices of the other partite class are labeled $\{b_1,\ldots,b_r\}$.
For $i \neq j$ color $(a_i,b_j)$ with the same color as the edge $(i,j)$ of $K_r$.
Color $(a_i,a_i)$ with any color incident with vertex $i$ in $K_r$.
The obtained coloring has the desired property.
Now, this construction together with the same ``blow up'' argument as in Lemma \ref{l31} yields
the following:
$$
f_k(K_{n,n}) \geq n(1-z_k)(1+o_n(1)).
$$
An upper bound density argument similar to the one in Theorem \ref{t3}
gives $f_k(K_{n,n}) \leq n(1-1/\sqrt{k})(1+o_n(1))$.
Recalling that $z_k=1/\sqrt{k}+O(1/k)$ we get that the upper and lower bounds for
$f_k(K_{n,n})$ are very close.

\section*{Acknowledgment}

The authors thank Yehuda Roditty for valuable discussions.

\end{document}